\documentclass[10pt]{article}
\usepackage{mathrsfs}
\usepackage{times}
\usepackage[T1]{fontenc}
\usepackage[dvips]{graphics}
\usepackage{epsfig}
\usepackage{amsmath,amsfonts,amsthm,amssymb,amscd}
\newtheorem{proposition}{Proposition}[section]

\newtheorem{theorem}{Theorem}[section]

\newtheorem{remark}[theorem]{Remark}

\def\phi{{\varphi}}

\DeclareSymbolFont{AMSb}{U}{msb}{m}{n}
\DeclareMathSymbol{\N}{\mathbin}{AMSb}{"4E}
\DeclareMathSymbol{\Z}{\mathbin}{AMSb}{"5A}
\DeclareMathSymbol{\R}{\mathbin}{AMSb}{"52}
\DeclareMathSymbol{\Q}{\mathbin}{AMSb}{"51}
\DeclareMathSymbol{\I}{\mathbin}{AMSb}{"49}
\DeclareMathSymbol{\C}{\mathbin}{AMSb}{"43}

\newcommand{\lam}{\lambda}

\newcommand{\alp}{\alpha}

\def\ni{\noindent}

\def\proof{{\ni \bf {Proof:} }}

\def\endproof{\hfill$\blacksquare$\vspace{6pt}}

\newcommand{\bsub}{\begin{subequations}}

\newcommand{\esub}{\end{subequations}$\!$}

\newtheorem{thm}{Theorem}[section]

\begin{document}

\title{Uniqueness of solutions for an elliptic equation 
modeling MEMS \textsc \bf }
\author{\sc{Pierpaolo ESPOSITO \footnote{Dipartimento di Matematica, Universit\`a degli Studi ``Roma Tre", 00146 Rome, Italy. E-mail: esposito@mat.uniroma3.it. Research supported by M.U.R.S.T., project ``Variational methods and nonlinear differential equations".}
\quad and 
\quad 
Nassif GHOUSSOUB\footnote{Department of Mathematics, University of British 
Columbia, Vancouver, B.C. Canada V6T 1Z2. E-mail: nassif@math.ubc.ca. 
Research partially supported by the Natural Science and Engineering 
Research Council of Canada.}}
\date{\today}
}
\maketitle

\section{Introduction}

We study the effect of the parameter $\lambda$, the dimension $N$, the 
profile $f$ and the geometry of the domain $\Omega \subset\mathbb{R}^N$, 
on the question of uniqueness  of the solutions  to the 
following elliptic boundary value problem with a singular nonlinearity:
$$\hskip 180 pt  \left\{ \begin{array}{ll} -\Delta
u= \frac{\lambda f(x)}{(1-u)^2} & \hbox{in }\Omega\\
0<u<1 &\hbox{in }\Omega\\
u=0 &\hbox{on }\partial \Omega. \end{array} \right. \hskip 130pt 
(S)_{\lambda, f} $$
This equation has been proposed as a model for a simple electrostatic 
Micro-Electromechanical System (MEMS) device consisting of a thin 
dielectric elastic membrane with boundary
supported at $0$ below a rigid ground plate located at height $z=1$. See 
\cite{P1,PB}.
A voltage -- directly proportional to the parameter $\lambda$ -- is 
applied, and the membrane deflects towards the ground plate and a 
snap-through may occur when
it exceeds a certain critical value $\lambda^*$, the pull-in voltage.

\medskip \noindent In \cite{JoLu} a fine ODE analysis of the radially 
symmetric case with a profile $f\equiv 1$ on a ball $B$, yields the 
following bifurcation diagram that describes the $L^\infty$-norm of the 
solutions $u$ -- which in this case necessarily coincides with $u(0)$ -- 
in terms of the corresponding voltage $\lambda$. \\

\begin{figure}[htbp]\label{B:fig2}
\begin{center}
{\includegraphics[width =
11cm,height=5.5cm,clip]{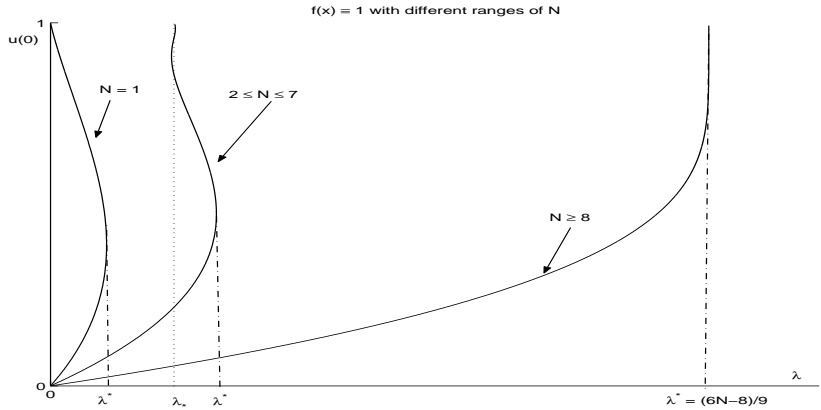}}
\caption{{\em   Plots of $u(0)$ versus $\lam$ for profile $f(x)
\equiv 1$
defined in the unit ball
$B_1(0)\subset \R^N$ with different ranges of $N$. In the case $N\ge
8$, we
have $\lam ^*=2(3N-4)/9$.} }
\end{center}
\end{figure}

\medskip \noindent The question whether the diagram above describes 
realistically the set of all solutions in more general domains and for 
non-constant profiles, and whether rigorous mathematical proofs can be 
given for such a description, has been the subject of many recent 
investigations. See  \cite{E1,EGG,GG1,GWA,GWB}.

\medskip \noindent We summarize in the following two theorems some of the 
established results concerning the above diagram. First, for every 
solution $u$ of $(S)_{\lambda, f}$, we consider the linearized operator   
$$
L_{u, \lam }=-\Delta -\frac{2\lam f}{(1-u)^{3}}
$$
and its eigenvalues $\{\mu_{k, \lambda}(u); k=1,2,\dots \}$ (with the 
convention that eigenvalues are
repeated according to their multiplicities). The Morse index $m(u, 
\lambda)$ of a solution $u$ is the largest $k$ for which 
$\mu_{k,\lambda}(u)$ is negative. A solution $u$ of $(S)_{\lambda, f}$ is 
said to be {\it stable} (resp., {\it semi-stable}) if $\mu_{1,\lam}( u 
)>0$ (resp., $\mu_{1, \lam}( u )\geq 0$).

\medskip \noindent A description of the first stable branch and of the 
higher unstable ones is given in the following.

\medskip \noindent {\bf Theorem A \cite{E1,EGG,GG1}} {\em Suppose $f$ is a 
smooth nonnegative
function in $\Omega$. Then, there exists a
finite $\lam ^*>0$ such that
\begin{enumerate}
\item If $0\leq \lam <\lam ^*$, there exists a (unique) minimal solution 
$u_\lambda$ of
$(S)_{\lambda, f}$ such that $\mu_{1,\lambda}(u_\lambda)>0$. It is also 
unique in the class of all semi-stable solutions.

\item If $\lam >\lam ^*$, there is no solution for $(S)_{\lambda, f}$.

\item If $1\leq N\leq 7$, then 
$u^*=\displaystyle\lim_{\lambda \uparrow \lambda^*}u_\lambda$ is a 
solution of $(S)_{\lambda^*, f}$ such that
$\mu_{1,\lambda^*}(u^*)=0,$
and  $u^*$  -- referred to as the extremal solution of problem 
$(S)_{\lambda,f}$ -- is the unique solution. 

\item If $1\leq N\leq 7$, there exists $\lambda_2^*$ with $0<\lambda_2^*<\lambda^*$ such that for any $\lambda \in (\lambda_2^*, 
\lambda^*)$, problem $(S)_{\lam,f }$ has a second solution $U_\lambda$ 
with $\mu_{1,\lambda}(U_\lambda) < 0$ and $\mu_{2,\lambda}(U_\lambda)> 0.$ 
Moreover, at $\lambda=\lambda_2^*$ there exists a second
solution  $U^*:=\displaystyle\lim_{\lambda \downarrow
\lambda_2^*} U_\lambda$ with
$$\mu_{1,\lambda_2^*}(U^*)<0  \quad {\rm and}\quad 
\mu_{2,\lambda_2^*}(U^*)=0.$$
\item Given a more specific potential $f$ in the form
\begin{equation}
\label{assf0}  f(x)=\left(\prod_{i=1}^k 
|x-p_i|^{\alpha_i}\right)h(x)\:,\quad \inf_\Omega h>0, 
\end{equation}
with points $p_i \in \Omega$, $\alpha_i \geq 0$, and given $u_n$ a 
solution of
$(S)_{\lam_n,f}$, we have the equivalence
$$\|u_n\|_\infty \to 1 \quad \Longleftrightarrow \quad m(u_n,\lambda_n) 
\to +\infty $$
as $n \to +\infty$. 
\end{enumerate}}
\noindent It was also shown in \cite{EGG} that the profile $f$ can 
dramatically change the bifurcation diagram, and totally alter the 
critical dimensions for compactness. Indeed, the following theorem 
summarizes the result related to the effect of power law profiles.

\medskip \noindent {\bf Theorem B \cite{EGG}} {\em Assume $\Omega$ is the 
unit ball $B$ and $f$ in the form
$$
f(x)=|x|^{\alpha} h(|x|)\:,\quad \inf_B h>0.
$$ 
Then we have

\begin{enumerate}

\item If $N \geq 8$ and 
$\alpha>\alpha_N:=\frac{3N-14-4\sqrt{6}}{4+2\sqrt{6}}$, the extremal 
solution $u^*$ is again a classical solution of $(S)_{\lam^*,f}$ such that 
$\mu_{1,\lambda^*}(u^*)=0$. 

\item If $N \geq 8$ and 
$\alpha>\alpha_N:=\frac{3N-14-4\sqrt{6}}{4+2\sqrt{6}}$, the conclusion of 
Theorem A-(4) still holds true.

\item On the other hand, if either $2\leq N\leq 7$ or $N\geq 8$, $0\leq 
\alpha \leq \alpha_N=\frac{3N-14-4\sqrt{6}}{4+2\sqrt{6}}$, for 
$f(x)=|x|^\alpha$ necessarily we have that 
$$u^*(x)=1-|x|^{\frac{2+\alpha}{3}}\,,\qquad 
\lambda^*=\frac{(2+\alpha)(3N+\alpha-4)}{9}.$$ 

\end{enumerate}}

\noindent The bifurcation diagram suggests the following conjectures:
\begin{enumerate}

\item For $2\leq N \leq 7$ there exists a curve $(\lambda(t),u(t))_{t\geq 
0}$ in the solution set 

\begin{equation} \label{solset}
\mathcal V=\Big\{(\lam,u)\in (0,+\infty) \times C^1(\bar \Omega):\: u 
\hbox{ is a solution of } (S)_{\lam,f } \Big\},
\end{equation}
starting from $(0,0)$ at $t=0$ and going to ``infinity": $\|u(t)\|_\infty 
\to 1$ as $t \to +\infty$, with infinitely many bifurcation or turning 
points in $\mathcal V$.

\item In dimension $N\geq 2$ and for any profile $f$, there exists a 
unique solution for small voltages $\lambda$. 

\item For $2\leq N \leq 7$ there exist exactly two solutions for $\lambda$ 
in a small left neighborhhod of $\lambda^*$.

\end{enumerate}

\medskip \noindent Conjectures $1$ and $2$ have been established for power 
law profiles in the radially symmetric case \cite{GWA}, and for the case 
where $f\equiv 1$, and $\Omega$ is a suitably symmetric domain in $\R^2$ 
\cite{GWB}. Indeed, in these cases Guo and Wei first show that
$$\lambda_*=\inf \{\lambda>0: \:(S)_{\lambda,f} \hbox{ has a non-minimal 
solution} \}>0,
$$
and then apply the fine bifurcation theory developed by Buffoni, Dancer 
and Toland \cite{BDT} to verify the validity of Conjecture $1$ too. 
Property $\lambda_*>0$ allows them to carry out some limiting argument and 
to prove that the Morse index of $u(t)$ blows up as $t \to +\infty$, which is 
crucial to show that infinitely many bifurcation or turning points occur 
along the curve. Thanks to Theorem A-(5), we shall be able in Section 2 to show  the validity 
of Conjecture $1$ in general domains $\Omega$,  by circumventing the need to prove that  $\lambda_*>0$. On the other hand, we shall prove in Section $3$ that   
indeed $\lambda_*>0$ for a large class of domains, and therefore we have uniqueness for small voltage. Our proofs simplify considerably those of Guo and Wei, and 
extend them  to general star-shaped domains $\Omega$ and power law profiles 
$f(x)=|x|^\alpha$, $\alpha \geq 0$.

\medskip \noindent Conjecture $3$ has been shown in \cite{E1} in the class 
of solutions $u$ with $m(u,\lambda)\leq k$, for every given $k\in \N$, and 
is still open in general.

\section{A quenching branch of solutions}

The first global result on the set of solutions in general domains was proved by the first 
author in \cite{E1}. By using a degree argument (repeated below),  he 
showed the following result. 

\begin{thm} Assume $2\leq N \leq 7$ and $f$ be as in (\ref{assf0}). There 
exist a sequence $\{\lambda_n \}_{n \in \N}$ and associated solution $u_n$ 
of $(S)_{\lam_n,f}$ so that 
$$m(u_n,\lam_n) \to +\infty \qquad \hbox{as } n \to +\infty.$$
\end{thm}

\noindent Let us introduce some notations according to Section 2.1 in 
\cite{BDT}. Set
$$X=Y=\{u \in C^1(\bar \Omega):\:u=0 \hbox{ on }\partial \Omega\}\:,\quad 
U=(0,+\infty)\times \{u \in X:\|u\|_\infty<1\},$$
and define the real analytic function $F: \R \times U \to Y$ as 
$F(\lambda,u)=u-\lambda K(u)$, where $K(u)=-\Delta^{-1} 
\left(f(x)(1-u)^{-2} \right)$ is a compact operator on every closed subset 
in $\{u \in X:\|u\|_\infty<1\}$ and
$\Delta^{-1}$ is the Laplacian resolvent with homogeneous Dirichlet
boundary condition. The solution set $\mathcal V$ given in (\ref{solset}) 
rewrites as
$$\mathcal V=\{(\lambda,u) \in U:\: F(\lambda,u)=0\},$$
and the projection of $\mathcal V$ onto $X$ is defined as
$$\Pi_X \mathcal V=\{u \in X:\: \exists \:\lam \hbox{ so that } (\lam ,u) 
\in \mathcal V\}.$$

\proof In view of Theorem A-(5), we have the equivalence
$$\sup_{(\lam,u) \in \mathcal V} \displaystyle \max_\Omega u=1
\qquad \Longleftrightarrow \qquad \sup_{(\lam,u) \in \mathcal V} 
m(u,\lambda)=+\infty.$$
Arguing by contradiction, we can assume that

\begin{equation} \label{compV}\sup_{(\lam,u) \in \mathcal V} \displaystyle 
\max_\Omega u\leq 1-2\delta,\qquad \sup_{(\lam,u) \in \mathcal V} 
m(u,\lambda)<+\infty 
\end{equation}
for some $\delta\in (0,\frac{1}{2})$. By Theorem 1.3 in \cite{E1} one can 
find $\lam_1,\lam_2 \in (0,\lam^*)$, $\lam_1<\lam_2$, so that 
$(S)_{\lam,f}$ possesses

\begin{itemize}

\item for $\lam_1$, only the (non degenerate) minimal solution
$u_{\lam_1}$ which satisfies $m(u_{\lam_1},\lam_1)=0$;

\item  for $\lam_2$, only the two (non degenerate) solutions
$u_{\lam_2}$, $U_{\lam_2}$ satisfying $m(u_{\lam_2},\lam_2)=0$ and 
$m(U_{\lam_2},\lam_2)=1$, respectively.

\end{itemize}

Consider a $\delta$-neighborhood $\mathcal V_\delta$ of $\Pi_X \mathcal V$:
$$\mathcal V_\delta:=\{ u \in X: \: \hbox{dist}_X (u,\Pi_X \mathcal V)\leq 
\delta \}.
$$
Note that (\ref{compV}) gives that $\mathcal V$ is contained in a closed 
subset of $\{u \in X:\|u\|_\infty<1\}$:
$$\mathcal V_\delta \subset \{u \in X:\:\|u\|_\infty \leq 1-\delta \}.
$$
We can now define the Leray-Schauder degree $d_\lam$ of $F(\lam,\cdot)$ on 
$\mathcal V_\delta$ with respect to zero, since by definition of $\Pi_X 
\mathcal V$ (the set of all solutions) $\partial
\mathcal V_\delta$ does not contain any solution of $(S)_{\lam,f}$ for
any value of $\lam$. Since $d_\lam$ is well defined for any $\lam
\in [0,\lam^*]$, by homotopy $d_{\lam_1}=d_{\lam_2}$.
\medskip \noindent To get a contradiction, let us now compute $d_{\lam_1}$ 
and
$d_{\lam_2}$. Since the only zero of $F(\lam_1,\cdot)$ in $\mathcal 
V_\delta$ is $u_{\lam_1}$ with Morse index zero, we have $d_{\lam_1}=1$. 
Since $F(\lam_2,\cdot)$ has in $\mathcal V_\delta$
exactly two zeroes $u_{\lam_2}$ and $U_{\lam_2}$ with Morse index
zero and one, respectively, we have $d_{\lam_2}=1-1=0$. This
contradicts $d_{\lam_1}=d_{\lam_2}$, and the proof is complete.
\endproof

\noindent We can now combine Theorem A-(5) with the fine bifurcation 
theory in \cite{BDT} to establish a more precise multiplicity result. See 
also \cite{Da}.\\ 
Observe that $\mathcal A_0 :=\{(\lambda,u_\lambda): \: \lambda \in 
(0,\lambda^*) \}$
is a maximal arc-connected subset of
$$S:=\{(\lambda,u)\in U:\:F(\lambda,u)=0 \hbox{ and }\partial_u F(\lambda, 
u):X \to Y \hbox{ is invertible} \}$$
with $\mathcal A_0 \subset S$. Assume that the extremal solution $u^*$ is 
a classical solution so to have $u^* \in (\bar S \cap U)\setminus S$. 
Assumption (C1) of Section 2.1 in \cite{BDT} does hold in our case. As far 
as condition (C2):
$$\{(\lambda,u)\in U:F(\lambda,u)=0 \} \hbox{ is open in }\{(\lambda,u) 
\in \R \times X:F(\lambda,u)=0 \},$$
let us stress that it is a weaker statement than requiring $U$ to be an 
open subset in $\R \times X$. In our case, the map $F(\lambda,u)$ is 
defined only in $U$ (and not in the whole $X$), and then condition (C2) 
does not make sense. However, we can replace it with the new condition 
(C2):
$$U \hbox{ is an open set in }\R \times X,$$
which does hold in our context. Since (C2) is used only in Theorem 
2.3-(iii) in \cite{BDT} to show that $S$ is open in $\bar S$, our new 
condtion (C2) does not cause any trouble in the arguments of \cite{BDT}.\\
Since $\partial_u F(\lambda,u)$ is a Fredholm operator of index $0$, by a 
Lyapunov-Schmidt reduction we have that assumptions (C3)-(C5) do hold in 
our case (let us stress that these conditions are local and $U$ is an open 
set in $\R \times X$).\\ 
Setting $\bar \lambda=0$ and defining the map $\nu: U\to [0,+\infty)$ as 
$\nu(\lambda,u)=\frac{1}{1-\|u\|_\infty}$, conditions (C6)-(C8) do hold in 
view of the property $\lambda \in [0,\lambda^*]$. Theorem 2.4 in 
\cite{BDT} then applies and gives the following.

\begin{thm} \label{ggh} Assume $u^*$ a classical solution of 
$(S)_{\lambda^*,f}$. Then there exists an analytic curve $(\hat 
\lambda(t),\hat u(t))_{t\geq 0}$ in $\mathcal V$ starting from $(0,0)$ and 
so that $\|\hat u(t)\|_\infty \to 1$ as $t \to +\infty$. Moreover, $\hat 
u(t)$ is a non-degenerate solution of $(S)_{\hat \lambda(t),f}$ except at 
isolated points.    

\end{thm} 

\noindent By the Implicit Function Theorem, the curve $(\hat \lambda(t), 
\hat u(t))$ can only have isolated intersections. If we now use the usual 
trick of finding a minimal continuum in $\{( \hat \lambda(t),\hat u(t)): 
\: t \geq 0\}$ joining $(0,0)$ to ``infinity", we obtain a continuous 
curve $(\lambda(t),u(t))$ in $\mathcal V$ with no self-intersections which 
is only piecewise analytic. Clearly, $\partial_u F(\lambda,u):X \to Y$ is 
still invertible along the curve except at isolated points.

\medskip \noindent Let now $2\leq N \leq 7$ and $f$ be as in 
(\ref{assf0}). By the equivalence in Theorem A-(5) we get that 
$m(\lambda(t),u(t)) \to +\infty$ as $t \to +\infty$, and then 
$\mu_{k,\lambda(t)}(u(t))<0$ for $t$ large, for every $k \geq 1$. Since 
$\mu_{k,\lambda(0)}(u(0))=\mu_{k,0}(0)>0$ and $u(t)$ is a non-degenerate 
solution of $(S)_{\lambda(t),f}$ except at isolated points, we find 
$t_k>0$ so that $\mu_{k,\lambda(t)}(u(t))$ changes from positive to 
negative sign across $t_k$. Since $\mu_{k+1,\lambda(t)}(u(t)) \geq 
\mu_{k,\lambda(t)}(u(t))$, we can choose $t_k$ to be non-increasing in $k$ 
and to have $t_k \to +\infty$ as $k \to +\infty$.\\ 
To study secondary bifurcations, we will use the gradient structure in the 
problem. Setting $(\lambda_k,u_k):=(\lambda(t_k),u(t_k))$, we have that 
$(\lambda_k, u_k) \notin S$. Choose $\delta>0$ small so that 
$\|u_k\|_\infty<1-\delta$, and replace the nonlinearity $(1-u)^{-2}$ with 
a regularized one:
$$f_\delta(u)=\left\{ \begin{array}{ll}(1-u)^{-2}&\hbox{if }u\leq 1-\delta
\,,\\[3mm]
\delta^{-2}&\hbox{if } u \geq 1-\delta,\end{array}\right.$$
and the map $F(\lambda,u)$ with the corresponding one 
$F_\delta(\lambda,u)$. We replace $X$ and $Y$ with $H^2(\Omega) \cap 
H_0^1(\Omega)$ and $L^2(\Omega)$, respectively. The map 
$F_\delta(\lambda,u)$ can be considered as a map from $\R \times X \to Y$ 
with a gradient structure:
$$\partial_u \mathcal J_\delta(\lambda,u)[\phi]=\left\langle 
F_\delta(\lambda,u),\phi \right\rangle_{L^2(\Omega)}$$
for every $\lambda \in \R$ and $u,\phi \in X$, where $\mathcal J_\delta:\R 
\times X \to \R$ is the functional given by
$$\mathcal J_\delta(\lambda,u)=\frac{1}{2}\int_\Omega |\nabla 
u|^2\:dx-\lambda \int_\Omega f(x) G_\delta(u)\:dx \:,\qquad 
G_\delta(u)=\int_0^u f_\delta(s)ds.$$
Assumptions (G1)-(G2) in Section 2.2 of \cite{BDT} do hold. We have that 
$(\lambda(t),u(t)) \in S$ for $t$ close to $t_k$ and $m(\lambda(t),u(t))$ 
changes across $t_k$. If $\lambda(t)$ is injective, by Proposition 2.7 
in \cite{BDT} we have that $(\lambda(t_k),u(t_k))$ is a bifurcation point. 
Then we get the validity of Conjecture $1$ as claimed below.

\begin{thm} \label{ecchime} Assume $2\leq N \leq 7$ and $f$ be as in 
(\ref{assf0}). Then there exists a continuous, piecewise analytic curve 
$(\lambda(t),u(t))_{t\geq 0}$ in $\mathcal V$, starting from $(0,0)$ and 
so that $\|\hat u(t)\|_\infty \to 1$ as $t \to +\infty$, which has either 
infinitely many turning points, i.e. points
where $(\lambda(t),u(t))$ changes direction (the branch locally ``bends 
back"), or infinitely many bifurcation points.

\end{thm}

\begin{remark} \rm In \cite{GWA} the above analysis is performed in the 
radial setting to obtain a curve $(\lambda(t),u(t))_{t\geq 0}$, as given 
by Theorem \ref{ecchime}, composed by radial solutions and so that 
$m_r(\lambda(t),u(t)) \to +\infty$ as $t \to +\infty$, $m_r(\lambda,u)$ 
being the radial Morse index of a solution $(\lambda,u)$. In this way, it 
can be shown that bifurcation points can't occur and then 
$(\lambda(t),u(t))_{t\geq 0}$ exhibits infinitely many turning points. 
Moreover, they can also deal with the case where $N\geq 8$ and 
$\alpha>\alpha_N$.
\end{remark}

\section{Uniqueness of solutions for small voltage in 
star-shaped domains}

We address the issue of uniqueness of solutions of the 
singular elliptic problem
\begin{equation}\label{p:uni:1}
\arraycolsep=1.5pt \left\{ \begin{array}{ll} -\Delta
u=\frac{\lambda|x|^\alp}{(1-u)^2}  & \hbox{in } \Omega \\
0<u<1 &\hbox{in }\Omega \\
u=0 &\hbox{on } \partial \Omega , \quad
\end{array} \right. \end{equation}
for $\lam>0$ small, where $\alp \ge 0$ and $\Omega$ is a bounded domain in  
$\mathbb{R}^N$, $N\geq 2$.
We shall make crucial use of the following extension of Pohozaev's 
identity due to Pucci and Serrin \cite{PS}.

\begin{proposition} Let $v$ be a solution of the boundary value problem 
$$\arraycolsep=1.5pt 
\left\{ \begin{array}{ll} -\Delta v=f(x,v) &\hbox{in } \Omega\\
v=0 &\hbox{on } \partial \Omega.
\end{array} \right. $$
Then for any $a\in \R$ and any $h\in C^2(\Omega; \R^N)\cap C^1(\bar 
\Omega; \R^N) $, the following identity holds
\begin{eqnarray}\label{poho}
\int_\Omega \left[ {\rm div}(h)F(x,v)-a v f(x,v)+\langle \nabla_x F(x,v), 
h \rangle \right]\, dx &=&\int_\Omega \left[(\frac{1}{2}{\rm 
div}(h)-a)|\nabla v|^2-\langle Dh \nabla v, \nabla v \rangle \right]\, dx 
\nonumber\\
&& + \frac{1}{2}  \int _{\partial\Omega}
|\nabla v|^2 \langle h,\nu \rangle d\sigma,
\end{eqnarray}
where $F(x,s)=\int_0^s f(x,t)\, dt$.
\end{proposition}
An application of the method in \cite{Sc} leads to the following result.  

\begin{thm} \label{p:uni:thm1} Let $\Omega \subset\mathbb{R}^N$ be a  
star-shaped domain
with respect to $0$. If $N\ge 3$, then for $\lambda$ small (\ref{p:uni:1}) 
has the unique solution $u_\lam$.
\end{thm}

\proof  Since $u_\lambda$ is the minimal solution of (\ref{p:uni:1}) for 
$\lam \in (0,\lam ^*)$, setting $v=u-u_\lam$ equation (\ref{p:uni:1}) 
rewrites equivalently as

\begin{equation}
\label{p:uni:v}\arraycolsep=1.5pt \left\{ \begin{array}{ll} -\Delta
v=\lambda |x|^\alp g_\lambda (x,v) &\hbox{in } \Omega \\
0\le v<1-u_\lam &\hbox{in } \Omega \\
v=0 &\hbox{on } \partial \Omega,  \quad
\end{array} \right. \end{equation}
where
\begin{equation}\label{p:uni:100}
g_\lambda(x,s)=\frac{1}{(1-u_\lam(x) -s)^2}-\frac{1}{(1-u_\lam(x) )^2}\,.
\end{equation}
It then  suffices to prove that the solutions of (\ref{p:uni:v}) must be 
trivial for $\lam $ small enough.
\medskip \noindent First compute $G_\lambda(x,s)$:
$$G_\lambda(x,s)=\int ^s_0 g_\lambda(x,t)\:dt=\frac{1}{1-u_\lam(x)
-s}-\frac{1}{1-u_\lam(x)}-\frac{s}{(1-u_\lam(x) )^2}.$$
Since the validity of the relation
$$\nabla_x \Big(|x|^\alp G_\lambda(x,s)\Big)=\alp |x|^{\alp -2}x 
G_\lambda(x,s)+|x|^\alp \nabla_x G_\lambda(x,s),$$
for $h(x)=\frac{x}{N}$ and $f(x,v)=|x|^\alpha g_\lambda(x,v)$ we apply the 
Pohozaev identity (\ref{poho}) to a solution $v$ of (\ref{p:uni:v}) to get
\begin{eqnarray}\label{p:uni:3}
&& \lam \displaystyle \int _\Omega  |x|^\alp 
\big[(1+\frac{\alpha}{N})G_\lambda(x,v(x))-av(x)g_\lambda(x, v(x))+
\langle \nabla_x G_\lambda (x,v(x)),\frac{x}{N} \rangle \big]\: dx 
\nonumber \\
&&=\displaystyle \int _\Omega
\big[(\frac{1}{2}-a )|\nabla
v|^2- \langle D(\frac{x}{N}) \nabla v,\nabla
v \rangle\big]dx+\displaystyle \frac{1}{2N}\int _{\partial\Omega }
|\nabla v|^2\langle x,\nu\rangle \: d\sigma \\
&&\ge (\frac{1}{2}-a- \frac{1}{N})\int _\Omega |\nabla v|^2dx.\nonumber
\end{eqnarray}
Since easy calculations show that 
$$\frac{G_\lambda(x,s)}{g_\lambda(x,s)}=\frac{1-u_\lam(x)
-s-\frac{(1-u_\lam(x) -s)^2(1-u_\lam(x) +s)}{(1-u_\lam(x)
)^2}}{1-\frac{(1-u_\lam(x) -s)^2}{(1-u_\lam(x) )^2}}$$
and 
$$\frac{\nabla_x 
G_\lambda(x,s)}{g_\lambda(x,s)}=\frac{1-\frac{(1-u_\lam(x)-s)^2(1-u_\lam(x) +2s)}{(1-u_\lam(x) )^3} }{1-\frac{(1-u_\lam(x) 
-s)^2}{(1-u_\lam(x) )^2}} \nabla u_\lam(x),$$
we obtain  
\begin{equation}
\label{p:uni:5}
\Big|\frac{G_\lambda(x,s)}{g_\lambda(x,s)}\Big
|\leq C_0|1-u_\lam 
(x)-s|\qquad 
\hbox{and}
\qquad
\Big|\frac{\nabla_x G_\lambda(x,s)}{g_\lambda(x,s)}-\nabla
u_\lam\Big|\leq C_0 |1-u_\lam (x)-s|^2|\nabla u_\lambda| 
\end{equation}
for some $C_0>0$, provided $\lambda$ is away from $\lambda^*$. Since 
$u_\lam \to 0$ in $C^1(\bar\Omega )$ as $\lam \to
0^+$, for $a>0$ from (\ref{p:uni:5}) we deduce that for any $(x,s)$ 
satisfying $|1-u_\lam (x) -s|\le \delta$

\begin{eqnarray}
&&(1+\frac{\alpha}{N}) 
G_\lambda(x,s)-as g_\lambda(x,s)
+\langle \nabla_x 
G_\lambda(x,s), \frac{x}{N}\rangle \label{p:uni:6}\\
&&\le  g_\lambda(x,s) \Big[C_0(1+\frac{\alpha}{N})\delta -a(1-u_\lam (x) 
-\delta )+\langle \nabla
u_\lam ,\frac{x}{N}\rangle +\frac{C_0}{N}\delta^2 |\nabla u_\lambda| 
|x|\Big] \leq 0, \nonumber\end{eqnarray}
provided $\delta$ and $\lambda$ are sufficiently small (depending on $a$). 
Since $N\geq 3$, we can pick $0<a <\frac{1}{2}-\frac{1}{N}$, and then by 
(\ref{p:uni:3}), (\ref{p:uni:6}) get that
\begin{eqnarray}
\label{p:uni:8}
&& \lam \int _{\{0\le v\le 1-u_\lam -\delta \}}
|x|^\alp \big[(1+\frac{\alp}{N})G_\lambda(x,v(x))-av(x)g_\lambda(x,v(x))+\langle \nabla_x 
G_\lambda(x,v(x)),\frac{x}{N}\rangle \big]dx\\
&&\ge (\frac{1}{2}-a-\frac{1}{N})\int _\Omega|\nabla
v|^2\: dx\ge C_s (\frac{1}{2}-a-\frac{1}{N}) \int_\Omega v^2\:dx
\nonumber 
\end{eqnarray}
for $\delta$ and $\lambda$ sufficiently small, where $C_s$ is the best 
constant in the Sobolev embedding of $H^1_0(\Omega)$ into $L^2(\Omega)$. 

\medskip \noindent On the other hand, since $G_\lambda(x,s)$, 
$sg_\lambda(x,s)$ and
$\nabla_x G_\lambda(x,s)$ are quadratic with respect to $s$ as $s\to 0$ 
(uniformly in $\lambda$ away from $\lambda^*$),
there exists a constant $C_\delta >0$ such that
\begin{equation}\label{cute}
(1+\frac{\alp}{N})G_\lambda(x,v(x))- avg_\lambda(x,v(x))+\langle \nabla_x 
G_\lambda(x,v(x)),\frac{x}{N}\rangle  \le C_\delta v^2(x)  
\end{equation}
for $x \in \{0\leq v\leq 1-u_\lambda-\delta \}$, uniformly for $\lambda$ 
away from $\lambda^*$. Combining (\ref{p:uni:8}) and (\ref{cute}) we get 
that
$$C_s\displaystyle
\big(\displaystyle\frac{1}{2}-a-\frac{1}{N}\big)\int _{\{0\le v\le
1-u_\lam -\delta \}} v^2dx\le \lam C_\delta \int _{\{0\le v\le
1-u_\lam -\delta \}}|x|^\alp v^2dx.$$ 
Therefore, for $\lambda$ sufficiently small we conclude that
$v\equiv 0$ in $\{0\le v \le 1-u_\lam -\delta \}$. This implies that
$v\equiv 0$ in $\Omega $ for sufficiently small $\lam$, and we
are done. \endproof\\

\noindent We now refine the above argument so as to cover other 
situations. To this aim, we consider  the -- potentially empty -- set
$$H(\Omega )=\Big\{h\in \mathcal{C}^1(\bar\Omega ,\mathbb{R}^N): {\rm 
div}(
h)\equiv 1 \hbox{ and $\langle h,\nu  \rangle \ge 0$ on $\partial 
\Omega$}\Big\},$$
and  the corresponding parameter
$$M(\Omega ):=\inf \Big\{\sup _{x\in \Omega}\bar\mu (h,x):\, h\in
H(\Omega )\Big\},$$
where
$$\bar \mu (h,x)=\frac{1}{2}\sup _{|\xi
|=1}\langle(Dh(x)+Dh(x)^T)\xi ,\xi
\rangle.$$
The following is an extension of Theorem \ref{p:uni:thm1}.

\begin{thm} \label{p:uni:thm3} Let $\Omega$ be a bounded domain in 
$\mathbb{R}^N$ such that $M(\Omega )<\frac{1}
{2}$. Then, for $\lambda$ small the minimal solution $u_\lam$ is the 
unique solution of problem (\ref{p:uni:1}), provided either $N\ge 3$ or 
$\alp >0$.\end{thm}

\proof As above, we shall prove that equation (\ref{p:uni:v}), with 
$g_\lambda$ as in (\ref{p:uni:100}), has only trivial solutions for 
$\lambda$ small. For a solution $v$ of (\ref{p:uni:v}) the Pohozaev 
identity (\ref{poho}) 
with $h\in H(\Omega )$ yields
\begin{eqnarray} \label{p:uni:3?}
&& \lam \int _\Omega  |x|^\alp \big[G_\lambda(x,v(x))(1+\alpha \langle 
\frac{x}{|x|^2}, h\rangle) -av(x)g_\lambda(x,v(x))+
\langle \nabla_x G_\lambda(x,v(x)), h\rangle \big]\:dx  \nonumber \\
&&= \int _\Omega
\big[(\frac{1}{2}-a)|\nabla
v|^2-\frac{1}{2} \langle (D h+D h^T)\nabla v,\nabla
v\rangle\big]dx+\displaystyle \frac{1}{2} \int _{\partial\Omega }
|\nabla v|^2\langle h,\nu\rangle\: d \sigma \\
&&\ge \int _\Omega (\frac{1}{2}-a-\bar
\mu (h,x)\big)|\nabla v|^2\:dx. \nonumber
\end{eqnarray}
Fix $0<a<\frac{1}{2}-M(\Omega )$ and choose
$h\in H(\Omega )$ such that
$$\frac{1}{2}-a-\sup _{x\in \Omega}\bar \mu (h,x)>0.$$
It follows from (\ref{p:uni:5}) that for any $(x,s)$ satisfying 
$|1-u_\lam(x) -s|\le \delta
|x|$ there holds
\begin{eqnarray}
&& G_\lambda(x,s)(1+\alpha \langle \frac{x}{|x|^2}, 
h\rangle)-avg_\lambda(x,s)+\langle \nabla_x G_\lambda(x,s),h\rangle 
\nonumber \\ 
&&\le g_\lambda(x,s)\big[C_0\delta |x|+\alpha C_0 \delta |h|-a(1-u_\lam 
-\delta |x|)+ \langle \nabla
u_\lam ,h\rangle
+C_0 \delta^2 |x|^2 |\nabla u_\lambda| |h| \big] \leq 0 \label{p:uni:1000}
\end{eqnarray}
provided $\lambda$ and $\delta$ are sufficiently small.
It then follows from (\ref{p:uni:3?}) and (\ref{p:uni:1000}) that
\begin{eqnarray}
\label{p:uni:11} 
&& \lam \int _{\{0\le v\le 1-u_\lam -\delta |x|\}}
|x|^\alp \Big[G_\lambda(x,v(x))(1+\alpha \langle \frac{x}{|x|^2}, 
h\rangle)- av(x)g_\lambda(x,v(x)) +\langle \nabla_x 
G_\lambda(x,v(x)),h\rangle \Big]\: dx \nonumber\\ 
&& \ge (\frac{1}{2}-a- \sup _{x\in \Omega}\bar \mu (h,x)) 
\int_\Omega|\nabla v|^2 \:dx.
\end{eqnarray}
On the other hand, there exists a constant $C_\delta >0$ such that
\begin{eqnarray*}
&&G_\lambda(x,v(x))(1+\alpha \langle \frac{x}{|x|^2}, 
h(x)\rangle)
-av(x)g_\lambda(x,v(x))
+<\nabla_x G_\lambda(x,v(x)),h(x)>\\
&&= \frac{v^2(x)}{(1-u_\lam(x) -v(x))(1-u_\lam(x)
)^2} 
(1+\alpha \langle \frac{x}{|x|^2}, h(x)\rangle) 
+\frac{av^2(x)[v(x)-2+2u_\lam(x))]}{(1-u_\lam(x)
-v(x))^2(1-u_\lam (x))^2}\\
&&+ \frac{v^2(x)(3-3u_\lambda(x)-2v(x))}{(1-u_\lam(x) 
-v(x))^2(1-u_\lam(x))^3}
<\nabla u_\lam(x) ,h(x)> \le C_\delta 
\frac{v^2(x)}{|x|^2}
\end{eqnarray*}
for $x \in \{0\leq v \leq 1-u_\lambda-\delta |x|\}$, uniformly for 
$\lambda$ away from $\lambda^*$.\\
If now $N\geq 3$, then Hardy's inequality combined with (\ref{p:uni:11}) 
implies
$$\frac{(N-2)^2}{4}(\frac{1}{2}-a-\sup
_{x\in \Omega} \bar \mu (h,x)\big)\int
_{\{0\le v\le 1-u_\lam -\delta |x|\}} \frac{v^2}{|x|^{2}}\: dx\le \lam 
C_\delta \int _{\{0\le
v\le 1-u_\lam -\delta |x|\}}\frac{v^2}{|x|^2}\:dx.
$$
On the other hand, when $N=2$ the space $H_0^1(\Omega)$ embeds continously 
into $L^p(\Omega)$ for every $p>1$, and then, by H\"older inequality, for 
$\alpha>0$ we get that
$$\int_\Omega \frac{v^2}{|x|^{2-\alpha}} \: dx \leq \left(\int_\Omega 
|x|^{-(2-\alpha)\frac{p}{p-2}}\:dx \right)^{\frac{p-2}{p}} 
\left(\int_\Omega |v|^{p} \:dx \right)^{\frac{2}{p}}\leq C_{N,\alpha}^{-1} 
\int_\Omega |\nabla v|^2\:dx
$$
provided $(2-\alpha)\frac{p}{p-2}<2$, which is true for $p$ large 
depending on $\alpha$ (see \cite{GM} for some very general Hardy 
inequalities). It combines with (\ref{p:uni:11}) to yield
$$C_{N,\alp }(\frac{1}{2}-a-\sup
_{x\in\Omega} \bar \mu (h,x))\int_{\{0\le v\le 1-u_\lam -\delta |x|\}}
\frac{v^2}{|x|^{2-\alp }}\:dx\le \lam C_\delta \int _{\{0\le
v\le 1-u_\lam -\delta |x|\}}\frac{v^2}{|x|^{2-\alp }}\: dx\,.
$$
In both cases, we can conclude that for $\lambda$ sufficiently small 
$v\equiv 0$ for $x \in \{0\le v\le 1-u_\lam -\delta |x|\}$, for some 
$\delta>0$ small. Since we can assume $\delta$ and $\lambda$ sufficiently 
small to have
$$1-u_\lam -\delta |x|\ge \frac{1}{2}\quad \mbox{in} \quad
\big\{x\in\Omega :\, |x|\ge \frac{1}{2}{\rm dist } (0,\partial
\Omega)\big\},$$
we then have $$ v\equiv 0 \quad \mbox{in} \quad \big\{x\in \Omega
:\, v(x)\le \frac{1}{2}\big\}\cap\big\{x\in \Omega :\, |x|\ge
\frac{1}{2}{\rm dist }(0,\partial \Omega)\big\}.
$$
Since $v=0$ on $\partial\Omega$ and the domain $\{x\in\Omega :\,
|x|\ge \frac{1}{2}{\rm dist }(0,\partial \Omega)\}$ is connected,  the
continuity of $v$ gives that
$$
v\equiv 0 \quad \mbox{in} \quad \big\{x\in \Omega :\, |x|\ge
\frac{1}{2}{\rm dist }(0,\partial \Omega)\big\}.
$$
Therefore, the maximum principle for elliptic equations implies
$v\equiv 0$ in $\Omega$, which completes the proof of Theorem
\ref{p:uni:thm3}.\endproof

\begin{remark} \rm In \cite{Sc} examples of dumbell shaped domains $\Omega 
\subset \R^N$ which satisfy condition $M(\Omega)<\frac{1}{2}$ are given 
for $N\geq 3$. When $N\geq 4$, there even exist topologically nontrivial 
domains with this property. Let us stress that in both cases $\Omega$ is 
not starlike, which means that the assumption $M(\Omega)<\frac{1}{2}$ on a domain $\Omega$
 is more general than being shar-shaped.
\end{remark}

\noindent The remaining case $N=2$ and $\alpha=0$, is a bit more delicate. We have the following result.
\begin{thm}  If $\Omega$ is either a strictly convex or a symmetric domain 
in $\R^2$,  then  $(S)_{\lambda, 1}$ has the unique solution $u_\lambda$ 
for small $\lam$. 
\end{thm}

\proof The crucial point here is the following inequality: for every 
solution $v$ of (\ref{p:uni:v}) there holds
$$\int_{\partial \Omega} |\nabla v|^2 \:d \sigma \geq l(\partial 
\Omega)^{-1} \left(\int_\Omega |\Delta v|\: dx \right)^2.$$
Indeed, we have that
$$\int_{\partial \Omega} |\nabla v|^2 \:d \sigma \geq l(\partial 
\Omega)^{-1} \left(\int_{\partial \Omega} |\nabla v|\:d \sigma \right)^2=
l(\partial \Omega)^{-1} \left(\int_{\partial \Omega} \partial_\nu v \: d 
\sigma \right)^2=
l(\partial \Omega)^{-1} \left(\int_\Omega |\Delta v|\:dx \right)^2,$$
where $l(\partial \Omega)$ is the length of $\partial \Omega$. Note that 
$-\Delta v=\lambda g_\lambda(x,v)\geq 0$ for every solution $u_\lambda+v$ 
of $(S)_{\lambda, 1}$, in view of the minimality of $u_\lambda$.

\medskip \noindent By Lemma 4 in \cite{Sc} for $\lambda$ small there 
exists $x_\lambda \in \Omega$ so that 

\begin{equation} \label{castoro}
\langle \nabla u_\lambda(x),x-x_\lambda \rangle\leq 0 \qquad \forall \:x 
\in \Omega.
\end{equation}
In particular, for $\lambda$ small $x_\lambda$ lies in a compact subset of 
$\Omega$ and, when $\Omega$ is symmetric, coincides exactly with the 
center of symmetries. In both situations, then we have that there exists 
$c_0>0$ so that
$$\langle x-x_\lambda,\nu(x)\rangle \geq c_0 \qquad \forall \: x \in 
\partial \Omega.$$
We use now the Pohozaev identity (\ref{poho}) with $a=0$ and 
$h(x)=\frac{x-x_\lambda}{2}$. For every solution $v$ of (\ref{p:uni:v}) it 
yields 
\begin{eqnarray}
\lam \int _\Omega  \big[G_\lambda(x,v(x))+
\langle \nabla_x G_\lambda (x,v(x)),\frac{x-x_\lambda}{2} \rangle \big]\: 
dx =\frac{1}{4}\int _{\partial\Omega }
|\nabla v|^2\langle x-x_\lambda,\nu\rangle \: d\sigma \ge \frac{c_0}{4} 
\left(\int_\Omega |\Delta v| \: dx \right)^2.\label{fft}
\end{eqnarray}
Since 
$$\nabla_x G_\lambda(x,s)=
(1-u_\lambda(x)-s)^{-2} \left[1-\frac{(1-u_\lam(x)
-s)^2(1-u_\lam(x) +2s)}{(1-u_\lam(x) )^3} \right] \nabla u_\lam(x),$$
by (\ref{castoro}) we easily see that
$$\langle \nabla_x G_\lambda(x,s), x-x_\lambda \rangle\leq 0$$
for $\lambda$ and $\delta$ small, provided $(x,s)$ satisfies 
$|1-u_\lambda(x)-s|\leq \delta$. Since $G_\lambda(x,s)$, $\nabla_x 
G_\lambda(x,s)$ are quadratic with respect to $s$ as $s\to 0$ (uniformly 
in $\lambda$ small),
there exists a constant $C_\delta >0$ such that
$$G_\lambda(x,v(x)) \le C_\delta v^2(x) \:,\qquad \langle \nabla_x 
G_\lambda(x,v(x)),\frac{x-x_\lambda}{2}\rangle  \le C_\delta v^2(x)$$ 
for $x \in \{0\leq v\leq 1-u_\lambda-\delta \}$, uniformly for $\lambda$ 
small.

Since on two-dimensional domains
$$\left(\int_\Omega |v|^p \: dx \right)^{\frac{1}{p}}\leq C_p \int_\Omega 
|\Delta v|\:dx$$
for every $p\geq 1$ and $v \in W^{2,1}(\Omega)$ so that $v=0$ on $\partial 
\Omega$, we get that
\begin{equation} \label{prima}
\lam \int _\Omega  \langle \nabla_x G_\lambda 
(x,v(x)),\frac{x-x_\lambda}{2} \rangle  
\: dx \leq \lambda C_\delta \int_\Omega v^2 \:dx\leq  \lambda C_\delta 
C_2^2 \left(\int_\Omega |\Delta v|\:dx\right)^2.
\end{equation}
As far as the term with $G_\lambda(x,v(x))$, fix $b \in (0,1)$ and split 
$\Omega$ as the disjoint union of $\Omega_1=\{v \leq b\}$ and 
$\Omega_2=\{v>b\}$. On $\Omega_1$ we have that
$$\lam \int _{\Omega_1}  G_\lambda(x,v(x))\: dx \leq \lambda C_\delta 
\int_\Omega v^2 \:dx \leq \lambda C_\delta C_2^2 \left(\int_\Omega |\Delta 
v|\:dx \right)^2$$
provided $\lambda$ and $\delta$ are small to satisfy $b \leq 
1-u_\lambda-\delta$ in $\Omega_1$.\\
Since for $\lambda$ small
$$\frac{G_\lambda(x,s)^2}{g_\lambda(x,s)} \leq C \quad \forall\: b\leq s 
\leq 1,$$
we have that
\begin{eqnarray*}
\lam \int _{\Omega_2}  G_\lambda(x,v(x))\: dx &\leq&  \lambda D_1 
\int_\Omega |v(x)|^{\frac{3}{2}} g_\lambda^{\frac{1}{2}}(x,v(x)) \:dx 
\leq  \lambda D_2 \left(\int_\Omega |v|^3 \: dx\right)^{\frac{1}{2}} 
\left(\int_\Omega g_\lambda(x,v(x))\:dx \right)^{\frac{1}{2}}\nonumber \\ 
&\leq&  \lambda^{\frac{1}{2}} D_3  \left(\int_\Omega |\Delta v|\:dx 
\right)^2
\end{eqnarray*}
for some positive constants $D_1$, $D_2$ and $D_3$. So we get that
\begin{eqnarray}\label{seconda}
\lam \int _\Omega  G_\lambda(x,v(x))\: dx \leq \left(\lambda C_\delta 
C_2^2+ \lambda^{\frac{1}{2}}D_3 \right) \left(\int_\Omega |\Delta v|\:dx 
\right)^2.
\end{eqnarray}
Inserting (\ref{prima})-(\ref{seconda}) into (\ref{fft}) finally we get
that
$$\left(2 \lambda C_\delta C_2^2+\lambda^{\frac{1}{2}}D_3  -\frac{c_0}{4} 
\right)\left(\int_\Omega |\Delta v|\:dx \right)^2 \geq 0,$$
and then $v \equiv 0$ for $\lambda$ small. \endproof

\end{document}